

\magnification=\magstep1
\overfullrule=0pt

\def\endpf{$\triangle$}
\baselineskip=15pt
\bigskip
\centerline{\bf On the Hilbert Schemes of Canonically-Embedded}
\centerline{\bf Curves of Genus 5 and 6}
\bigskip
\centerline{John B. Little}
\bigskip
\centerline{November 30, 1993}
\bigskip
\noindent
{\bf \S 1.  Introduction.}
\bigskip
Throughout this paper, the phrase {\it canonically-embedded curve of
genus $g$} will be used to refer to any pure 1-dimensional, non-degenerate
subscheme $C$ of ${\bf P}^{g-1}$ over an algebraically closed field
$k$, for which
$${\cal O}_C(1) \cong \omega_C, \hbox{  the dualizing sheaf}$$
and
$$h^0(C, {\cal O}_C) = 1, \qquad h^0(C, \omega_C) = g\ .$$
The singularities of $C$ (if any) are Gorenstein, and $C$ is connected
of degree $2g-2$ and arithmetic genus $g$.

In a recent paper ([S]), Schreyer has proved that Petri's normalization
of the homogeneous ideal $I(C)$ of a smooth canonically-embedded curve
(see also [P], [M], [SD], [ACGH] for the case of smooth curves) can be also
carried out for singular curves, provided that the curve has a simple
$(g-2)$-secant (i.e. a linear $W \cong {\bf P}^{g-3}$ intersecting $C$
transversely at
exactly $g-2$ (smooth) points).  We will call such $C$ {\it Petri-general}.
Moreover, he shows that the variety defined by any
system of quadrics and cubics of Petri's form is a Petri-general canonical
curve (Theorem 1.4 of [S]).  The Groebner basis techniques used in [S] also
lead naturally to a construction of a {\it Petri scheme} ${\cal P}_g$
parametrizing all ideals of Petri's form for a fixed $g$.

In the present paper (in effect an extended footnote to [S])
we will use the
map from ${\cal P}_g$ to the Hilbert scheme of curves of degree $2g-2$ and
arithmetic genus $g$ in ${\bf P}^{g-1}$ to study the low-genus cases $g = 5,6$.
Of course, the situation for smooth curves is completely understood for
all $g$.
Since the moduli space of smooth curves is irreducible, and the canonical
embedding is determined by a choice of basis in $H^0(C, \omega_C)$,
the points of the Hilbert scheme corresponding to smooth canonically-embedded
curves all lie on one irreducible component of dimension
$$\hbox{dim}({\cal M}_g) + \hbox{dim}(PGL(g)) = 3g - 3 + g^2 - 1\ .$$
(The Petri-general curves obtained as in [S] correspond to a certain
subscheme
of dimension $7g-7$, defined by incidence conditions with the fixed
$(g-2)$-secant used in Petri's construction.)
However, there can be other components of the Hilbert scheme whose
general points correspond to singular curves.  Our main results (see
Theorem 2.1, Theorem 3.1, Theorem 3.5) are
that this possibility {\it does not} arise (yet) for Petri-general
curves of $g=5$ or $g=6$:  If $C$ is Petri-general, then $[C]$ lies on the
component of the Hilbert scheme whose general point corresponds to a
smooth curve.

The $g = 5$ result, Theorem 2.1, is quite easy (and certainly not new,
though we could not find an explicit statement in the literature).  It
should follow, for example, from the algebraic structure theorems
for Gorenstein ideals of codimension 3 in [BE].  However, we will
use the extension of the Enriques-Petri theorem to singular curves
given in \S 3 of [S].
A canonical curve of genus 5 is either a complete intersection of
three quadrics in ${\bf P}^4$, or else it lies on
a surface of degree 3 in ${\bf P}^4$ (a rational normal cubic scroll, or a
degeneration)
cut out by the quadrics in $I(C)$.  Each of these
two possibilities yields an irreducible family of curves,
the ``non-trigonal'' complete intersections ${{\cal H}_5}'$ and the
``trigonal'' curves for which $I(C)$ is not generated by quadrics
${{\cal H}_5}''$, and ${{\cal H}_5}''$ is contained in the closure
of  ${{\cal H}_5}'$.

In the genus 6 case, the situation is very similar, but somewhat more
complicated.
The most difficult part of the proof, in fact, is to show that the family
of Petri-general canonical curves of genus 6, whose ideals are
generated by quadrics, is irreducible.  We do this by showing that each
of these curves is the complete intersection of a surface of degree
5 (a quintic Del Pezzo surface or a degeneration) and a quadric in ${\bf P}^5$.
The existence of such a surface for smooth $C$ is, of course, classical
(see the comments at the start of \S 3).  Our proof here uses many of
the ideas of Schreyer's study of the 2nd syzygy module of $I(C)$
in \S 4 of [S].

This paper was written while I was a visitor at ACSyAM, the Army
Center of Excellence for Symbolic Methods in Applied Mathematics, at the
Mathematical Sciences Institute at Cornell University.  I would to
thank Moss Sweedler for his hospitality, and Mike Stillman for several
very helpful conversations.  I would also like to acknowledge
Mira Bernstein, Diane Jamrog, Millie Niss, and Rachel Pries, the
members of an NSF REU group at the 1992 Regional Geometry Institute
at Amherst College who studied the genus 5 Petri scheme with me as
I made my first acquaintance with the results of [S].
\bigskip
\noindent
{\bf \S 2.  The Genus 5 Case.}
\bigskip
Let ${\cal H}_5$ be the open subvariety of the Hilbert scheme of
curves of degree 8 and arithmetic genus 5 in
${\bf P}^4$ consisting
of points corresponding to Petri-general curves $C$.  Recall, this
means that $C$ contains 3 smooth points in spanning a 2-plane
which intersects $C$ transversely exactly in those points.
(However, $C$ may have any other singularities, non-reduced
components, etc. consistent with this requirement.)
In this section we will prove that ${\cal H}_5$ is irreducible.
Although this fact is certainly well-known in a sense, we include it for
completeness, and for the way that some
of the ideas we use in this proof will reappear in the proofs
for genus 6 curves later in the paper.
\bigskip
\proclaim Theorem 2.1.
${\cal H}_5$ is irreducible of dimension 36, contained in one
irreducible component of the Hilbert scheme. \par
\bigskip
\noindent
{\bf Proof.}  The proof is in three sections.  First we deal with the
curves $C$ for which $I(C)$ is generated by quadrics (the ``non-trigonal''
case).  Since $\hbox{dim }I(C)_2 = 3$ for any genus 5 canonical curve,
we have that $C$ is a complete intersection, and the corresponding
open subvariety ${{\cal H}_5}'$ of ${\cal H}_5$ is isomorphic to a
Zariski-open dense subset
of the Grassmannian of three-dimensional subspaces of
$H^0({\bf P}^4, {\cal O}_{{\bf P}^4}(2))$.  Hence ${{\cal H}_5}'$ is
irreducible,
of dimension 36.

Second, we consider the cases for which $I(C)$ is not generated by quadrics
(the ``trigonal'' case).
Our main tool will be the extension of the Enriques-Petri Theorem to
Petri-general,
but possibly singular or reducible
curves in \S3 of [S].  Choosing coordinates \`a la Petri, we
have the Petri coefficient $\rho_{123} = 0$, and the quadrics in $I(C)$
have a basis of the form:
$$\eqalign{f_{12} &= x_1 x_2 - a_{112} x_1 - a_{212} x_2 - q_{12}(x_4, x_0)\cr
           f_{13} &= x_1 x_3 - a_{113} x_1 - a_{313} x_3 - q_{13}(x_4, x_0)\cr
           f_{23} &= x_2 x_3 - a_{223} x_2 - a_{323} x_3 - q_{23}(x_4, x_0)\cr}
\leqno(1)$$
Again since $\rho_{123} = 0$, the Petri syzygies on $I(C)$ (see Corollary 1.5
of [S])
read as follows:
$$x_j f_{ik} - x_k f_{ij} - \sum_{s=1,s\ne j}^{3} a_{sik}f_{sj}
 - \sum_{s=1, s\ne k}^{3} a_{sij} f_{sk} = 0.$$
It follows that $F = \{f_{12}, f_{13}, f_{23}\}$ forms a Groebner basis
for the ideal $J$ they generate (using the graded reverse lexicographic order,
with
the variables ordered $x_1 > x_2  > x_3 > x_4 > x_0$ as in [S]).

The variety $V(J)$ is a reduced, arithmetically
Cohen-Macaulay surface $S$ of degree 3 in ${\bf P}^4$.  (In the case that $C$
is a smooth curve, $S$ is always a smooth rational normal scroll.  However,
for singular $C$, in addition to smooth scrolls, cones over twisted cubic
curves and reducible surfaces can and do appear (see the examples in \S3 of
[S]).  We claim that the family ${\cal S}$ of all such surfaces in ${\bf P}^4$
is
irreducible.

One proof of this (well-known) fact follows directly from the Petri form (1)
for the generators of $J$.  Using a linear change of coordinates, we
can put the equations of any surface $S$ of the family
${\cal S}$ into this form.  The fact that $F$ forms a Groebner basis implies
that $S$ lies in a 1-parameter family flatly deforming to the surface
$S_0 = V(x_1 x_2, x_1 x_3, x_2 x_3)$.  The
Zariski tangent space to Hilb at $[S_0]$ has dimension equal to
$\hbox{dim}(H^0(N_{S_0 | {\bf P}^4}))$, which is equal to 18 by a direct
calculation.  This is exactly equal to the dimension
of the component of the Hilbert scheme of surfaces of degree 3 containing
the rational normal scrolls.  It follows that $[S_0]$ is a smooth
point of this Hilbert scheme, and hence that all of our surfaces
correspond to points on one irreducible component.

Using the results of \S 3 of [S], on each surface $S$,
the canonically-embedded ``trigonal'' curves form a
17-dimensional irreducible system.

By the discussion on p. 102 of [S], in addition to the
smooth rational normal scrolls, the other possible surfaces that
can occur here are:
\item{a)} $S = $ a cone over a twisted cubic curve,
\item{b)} $S = Q \cup P$, where $Q$ is a quadric in a ${\bf P}^3 \subset {\bf
P}^4$,
and $P$ is a ${\bf P}^2 \subset {\bf P}^4$ meeting $Q$ along a line (which
meets $C$ in three points), or
\item{c)} $S = P_1 \cup P_2 \cup P_3$, where $P_i$ are ${\bf P}^2$'s
intersecting
in the following way (after suitably numbering the components): $P_1$
and $P_3$ meet $P_2$ along lines, while $P_1 \cap P_3$ is a point.

Given a canonically-embedded ``trigonal'' curve $C$, the complete
intersection of the surface $S$ with one cubic not containing $S$
(e.g. one of the $G_{kl}$ in Petri's normalization)
will be $C$ union a line $L$ (a line of the ruling on a smooth scroll
or in case a), in $Q$ in case b), in $P_2$ in case c) above).
Conversely, if we fix such a line $L \subset S$, and consider the
linear system of cubics in ${\bf P}^4$ containing $L$ but
not containing $S$, the residual intersection
will be a canonically-embedded curve of genus 5 on $S$.  The dimension
of the space of cubics containing $L$, modulo cubics in $I(S)$
is:
$$({4 + 3\choose 3} - 4 - 13) - 1 = 17\ .$$

 The conclusion of this analysis is that the family ${{\cal H}_5}''$ of all
``trigonal'' canonically-embedded curves of genus 5 is irreducible of dimension
$18 + 17 = 35$.

Finally, we must show how ${{\cal H}_5}'$ and ${{\cal H}_5}''$ are related.  It
is easy to see that the smooth trigonal curves correspond to a dense
open subset of ${{\cal H}_5}''$, since the general surface $S$ in our family
${\cal S}$ is a smooth scroll, and the general element of the linear system
$|3b + 5f|$ on a smooth scroll is smooth.  Every smooth trigonal
curve is a limit of smooth non-trigonal curves.  Hence, an open dense
subset of the irreducible ${{\cal H}_5}''$ is contained in the Zariski closure
of
${{\cal H}_5}''$.  Hence all of ${{\cal H}_5}''$ is contained in the
closure of ${{\cal H}_5}'$, and this completes the proof of the proposition.
\endpf
\bigskip
\noindent
{\bf \S 3.  The Genus 6 Case.}
\bigskip
We now turn our attention to canonically-embedded curves $C$ of genus 6.
Recall that the
Enriques-Petri theorem for smooth curves of genus 6 says that $I(C)$ is
generated by
quadrics unless either $C$ is trigonal, or $C$ is isomorphic to a
smooth plane quintic.  In the first case, the quadrics cut out a rational
normal
scroll in ${\bf P}^5$ containing $C$; in the second case, the quadrics
cut out a Veronese surface.  In this section, we will begin by proving the
following
theorem about the Hilbert scheme of (possibly singular) canonically-embedded
curves
of genus 6.
\bigskip
\proclaim Theorem 3.1.  Let ${{\cal H}_6}'$ be the subset of the
Hilbert scheme of curves of arithmetic genus 6 and degree 10 in ${\bf P}^5$
parametrizing Petri-general curves $C$ for which $I(C)$ is generated by
quadrics.
Then ${\cal H'}_6$ is contained in one irreducible component of
Hilb. \par
\bigskip
Before beginning the proof (which is somewhat involved) we want to describe the
main idea and illustrate it with an example that led us to the general
statement.

 By a classical result (see [AH] for a modern treatment),
a smooth curve of genus 6, neither trigonal nor isomorphic
to a smooth plane quintic, lies on a
surface of degree 5 in ${\bf P}^5$ (a possibly degenerate quintic Del Pezzo
surface,
or a cone over a elliptic normal curve in ${\bf P}^4$).  The curve
is then the complete intersection of the surface and one further quadric.
For instance, projecting a {\it general} genus 6 canonical curve $C$ from a
general 4-secant 2-plane to $C$ (spanned by the points of one divisor
in one of the $g_4^1$'s on $C$) yields a plane curve of degree 6 with
four double points, birational
to $C$.  The canonical divisors on the plane model are cut by cubics passing
through the four double points.  This linear system of cubics maps
${\bf P}^2$ to a quintic Del Pezzo surface containing the original canonical
curve $C$.
This fact was also used in [KS] to show that the moduli space of
smooth curves with  $g = 6$ is stably rational.

What we will show is that even for singular curves, we still have a surface of
degree 5
containing $C$ that plays the same role.  To prove
Theorem 3.1 we will use the fact that the family of such surfaces is
irreducible.
Since the canonical curves on any one such surface are simply cut by the linear
system
of quadrics from ${\bf P}^5$, we will get the irreducibility statement to be
proved.

The main step of the proof will be to isolate the quadrics defining the
surface of degree 5.  The whole space of quadrics $I(C)_2$ is 6-dimensional,
and
as before we can take a basis in Petri form:
$$f_{ij} = x_i x_j - a_{1ij} x_1 - a_{2ij} x_2 - a_{3ij} x_3 - a_{4ij} x_4 -
q_{ij}(x_0, x_5),$$
where $1 \le i < j \le 4$, and if $i,j,k$ are distinct, $a_{kij} = \rho_{kij}
\alpha_k$ for
some constants $\rho_{kij}$, and non-zero
linear forms $\alpha_k = \alpha_k(x_0,x_5)$. We
will prove shortly that if $I(C)$ is generated by quadrics, then
at least two of the $\rho_{ijk}$ are different from zero.  By renumbering
the variables $x_1, x_2, x_3, x_4$ if necessary, we may assume $\rho_{123},
\rho_{124}$ are non-zero.  If this is the case, then
we will show that the (5-dimensional) subspace of $I(C)_2$ spanned by
$$\eqalign{F_1 &= \rho_{134} f_{12} - \rho_{123} f_{14}\cr
           F_2 &= \rho_{134} f_{23} - \rho_{123} f_{34}\cr
           F_3 &= \rho_{124} f_{13} - \rho_{123} f_{14}\cr
           F_4 &= \rho_{124} f_{23} - \rho_{123} f_{24}\cr
           F_5 &= \rho_{234} f_{12} - \rho_{123} f_{24}\cr
           F_6 &= \rho_{234} f_{13} - \rho_{123} f_{34}\cr}\leqno(2)$$
generates the ideal of a surface of degree 5 and sectional genus 1.
\bigskip\noindent
{\bf Remarks.}
Readers of [S] will note a clear
parallel between (2) and the quadrics in Claim 1 in the proof of Theorem 4.1
of that paper.  However, this initial $g = 6$ case seems to have some different
features from the $g \ge 7$ cases treated there, and we were not able
to deduce our desired result directly from the
techniques of \S4 of [S].  We should note that it is probably also possible to
derive this result from the purely algebraic structure theorem for
Gorenstein ideals of codimension 4 and deviation 2 in [HM],
but the proof using the Petri machinery is appealing in its own right, so we
proceed this way.
\bigskip
Here is an amusing singular example that gives one possible
type of singular surface that appears in this context.
\bigskip
\noindent
{\bf Example 3.2.}  Consider the ideal $I$ generated by the Petri-form
quadrics:
$$f_{ij} = x_i x_j - (x_k + x_l)(x_0 + x_5) - x_0 x_5,$$
where as usual $1 \le i < j \le 4$, and here $\{i,j,k,l\} = \{1,2,3,4\}$.
For this example, we can take $\rho_{ijk} = 1$, all $i,j,k$, and
$\alpha_k = x_0 + x_5$, all $k$.  Using a computer algebra system,
it is easy to see that $I$ has a Groebner
basis of the form described in Theorem 1.4 of [S], so $C = V(I)$ is
a canonical curve of arithmetic genus 6. It is not difficult
to see that $C$ is a union of five smooth conics $C_i$ in planes $P_i$ defined
by
$$\eqalign{P_1 &= V(x_2 + x_5 + x_0, x_3 + x_5 + x_0, x_4 + x_5 + x_0)\cr
           P_2 &= V(x_1 + x_5 + x_0, x_3 + x_5 + x_0, x_4 + x_5 + x_0)\cr
           P_3 &= V(x_1 + x_5 + x_0, x_2 + x_5 + x_0, x_4 + x_5 + x_0)\cr
           P_4 &= V(x_1 + x_5 + x_0, x_2 + x_5 + x_0, x_3 + x_5 + x_0)\cr
           P_5 &= V(x_1 - x_4, x_2 - x_4, x_3 - x_4)\cr}$$
The five planes all contain the line
$$L = V(x_1 + x_5 + x_0, x_2 + x_5 + x_0, x_3 + x_5 + x_0, x_4 + x_5 + x_0)$$
and the conics $C_i$ all meet $L$ in the same two points $p, q$.  (The tangents
to the $C_i$ at $p$ all lie in a hyperplane, so the singularity has
$\delta$-invariant 5; the situation at $q$ is the same.)  In this case the
surface of degree 5 defined by the combinations of the $f_{ij}$ given in
(2) above is $S = P_1 \cup \cdots \cup P_5$.  Note that the general hyperplane
section $S \cap H$ is a union of 5 concurrent lines spanning $H \cong
{\bf P}^4$,
a curve of arithmetic genus 1.
\bigskip
\noindent
{\bf Proof of Theorem 3.1.}  Let $C$ be a Petri-general canonically-embedded
curve of arithmetic genus 6 in ${\bf P}^5$, whose ideal is generated by
quadrics.
We will begin by proving the assertion above that at least two of the
Petri coefficients $\rho_{ijk}$ must be non-zero.  Let ${\cal T}$ be
the graph with vertices ${\cal V} = \{1,2,3,4\}$, and an edge
$(i,j)$ if and only
if there is some $k$ such that $\rho_{ijk} \ne 0$.  By Proposition 3.2
of [S], in the minimal free resolution of the homogeneous coordinate ring
of $C$, the graded Betti number $\beta_{13}$ (giving the number of
cubics in a minimal basis for $I(C)$) satisfies
$$\beta_{13} = \# \hbox{ connected components of } {\cal T} - 1$$
By our assumption, $\beta_{13} = 0$, so ${\cal T}$ must be connected.  By
the definition and the symmetry of the $\rho_{ijk}$ in the indices,
this implies that there is at most
one edge of the complete graph on ${\cal V}$ that is not contained in
${\cal T}$.  After
renumbering if necessary, the potentially missing edge can be
taken to be $(3,4)$, and hence $\rho_{123} \ne 0$ and $\rho_{124}
\ne 0$.  Under this assumption, the quadrics $F_i$ given in (2) above always
generate a 5-dimensional vector subspace of $I(C)_2$; there is exactly
one linear dependence between them:
$$\rho_{234}(F_1 - F_3) - \rho_{124}(F_2 - F_6) + \rho_{134}(F_4 - F_5) = 0\
.$$

\bigskip
\noindent
{\bf Step 1.}
\bigskip
Let $J = \langle F_1, \cdots, F_6\rangle$.  As explained above, the first
step in the proof will be to show that $V(J)$ is a surface of degree
5 in ${\bf P}^5$.  To do this, we will analyze the form of the unique
reduced Groebner basis
for $J$ with respect to the graded reverse lexicographic order with the
variables ordered $x_1 > x_2 > x_3 > x_4 > x_5 > x_0$.

The results depend on whether $\rho_{134}$ and $\rho_{234}$ are zero; we will
consider the case where all $\rho_{ijk} \ne 0$ first.  (Note that this
will be the case, for example, for a generic choice of $(g-2)$-secant in
Petri's
construction if $C$ is irreducible.)

For simplicity, we begin by taking linear combinations of the $F_i$ to
eliminate common terms and to isolate the possible leading monomials of
quadrics in the ideal.  The resulting basis for $J_2$ is:
$$\eqalign{F_1' &= \rho_{234}\rho_{134} f_{12} - \rho_{124} \rho_{123}
f_{34}\cr
           F_2' &= \rho_{234} f_{13} - \rho_{123} f_{34}\cr
           F_3' &= \rho_{234} f_{14} - \rho_{124} f_{34}\cr
           F_4' &= \rho_{134} f_{23} - \rho_{123} f_{34}\cr
           F_5' &= \rho_{134} f_{24} - \rho_{124} f_{34}\cr}$$
For instance, $F_1' = \rho_{234}(F_1 - F_3) + \rho_{124} F_6$.  We omit
the rest of the details in this calculation.  For future reference, however,
we note the following observation.
\bigskip
\proclaim Observation 3.3.  {\it All\/} quadrics of the forms
$$\rho_{ijk} f_{lj} - \rho_{ljk} f_{ij}$$
 and
$$\rho_{ikl} \rho_{jkl} f_{ij} - \rho_{ijk} \rho_{ijl} f_{kl},$$
(where $\{i,j,k,l\} = \{1,2,3,4\}$ in each case) belong to $J$.\par
\bigskip
This may be seen directly by forming linear combinations as above.
Now, applying Buchberger's algorithm, we begin the Groebner basis computation
on the $F_i'$.  At the first step, the $S$-polynomial $S(F_1', F_2')$ yields
$$x_3 F_1' - \rho_{134} x_2 F_2'
\equiv \rho_{123}(\rho_{234}x_2 x_3 x_4 - \rho_{124} x_3^2 x_4) \hbox{ mod }
\langle x_5, x_0 \rangle\ .$$
Replacing this last polynomial with its remainder on division by the
$F_i'$ and adjusting constants, we obtain a new Groebner basis element
$$G \equiv \rho_{123} x_3^2 x_4 - \rho_{124} x_3 x_4^2 \hbox{ mod }
\langle x_5, x_0 \rangle\ .$$

We claim that ${\cal G} = \{F_1', \cdots, F_5', G\}$ is the reduced Groebner
basis for $J$.  Indeed, working modulo $\langle x_5, x_0 \rangle$,
it is easy to see that all further $S$-pairs reduce to zero, modulo
$\langle x_5, x_0\rangle$.  Hence to prove
the claim, it suffices to show that $x_5, x_0$ are a regular sequence
in $k[x_1, \ldots, x_5, x_0]/J$, or that the resulting syzygies modulo
$\langle x_5, x_0 \rangle$ all lift to syzygies on ${\cal G}$. Using the
following lemma, we will show that this is a
consequence of the Petri syzygies on the generators of $I(C)$.
\bigskip
\proclaim Lemma 3.4.  Let $f_{ij}$ be a basis for the quadrics in the ideal
of a canonical curve of genus $g \ge 6$ in Petri's form, and let
$\{i,j,k,l\}$ be any four distinct indices in $\{1,2, \ldots, g-2\}$.  Then
each syzygy on the leading terms of the $f_{ij}$ of the form:
$$\eqalign{\rho_{iln} x_k (\rho_{ikn} x_i x_l - \rho_{ikl} x_i x_n) &+
  \rho_{ikn} x_l (\rho_{ikl} x_i x_n - \rho_{iln} x_i x_k) +\cr
  \rho_{ikl} x_n (\rho_{iln} x_i x_k &- \rho_{ikn} x_i x_l) = 0\cr}
\leqno(3)$$
lifts to a syzygy on the quadrics of the form
$$\rho_{\alpha\beta\gamma} f_{\delta\epsilon} - \rho_{\alpha\beta\delta}
f_{\gamma\epsilon}$$
in $I(C)$.\par
\bigskip
{\noindent}
{\bf Proof.}  We consider the Petri syzygies in the form
$$x_k f_{il} - x_l f_{ik} + \sum_{s \ne k} a_{sil} f_{sk} -
\sum_{s \ne l} a_{sik} f_{sl} + \rho_{ikl} G_{kl} = 0\leqno(S_{ikl})$$
where the $G_{kl}$ are Petri's cubics in $I(C)$, satisfying relations
$$G_{kl} + G_{ln} + G_{nk} = 0\ .\leqno(4)$$
Consider the linear combination
$$\rho_{ikn} \rho_{iln} S_{ikl} + \rho_{ikn} \rho_{ikl} S_{iln} +
\rho_{iln} \rho_{ikl} S_{ink}\leqno(5)$$
Using (4), the $G$ terms appearing in (5) cancel, so we obtain a relation on
quadrics only.  The remaining terms are:
$$\eqalign{0 =
\rho_{ikn}\rho_{iln}(x_k f_{il} - x_l f_{ik}&+
\sum_{s \ne k} a_{sil} f_{sk} - \sum_{s \ne l} a_{sik} f_{sl}) + \cr
\rho_{ikn}\rho_{ikl}(x_l f_{in} - x_n f_{il}&+
\sum_{s \ne l} a_{sin} f_{sl} - \sum_{s \ne n} a_{sil} f_{sn}) + \cr
\rho_{iln}\rho_{ikl}(x_n f_{ik} - x_k f_{in}&+
\sum_{s \ne n} a_{sik} f_{sn} - \sum_{s \ne k} a_{sin}f_{sl})\cr}$$
Rearranging and using the relations $a_{ijk} = \rho_{ijk}\alpha_i$
when $i,j,k$ are distinct, we have:
$$
\eqalign{&\rho_{iln} x_k (\rho_{ikn} f_{il} - \rho_{ikl} f_{in})
        +\rho_{ikn} x_l (\rho_{ikl} f_{in} - \rho_{iln} f_{ik})
        +\rho_{ikl} x_n (\rho_{iln} f_{ik} - \rho_{ikn} f_{il})\cr
= &\ \rho_{ikn} \sum_{s\ne k,n} a_{sil}(\rho_{ikl} f_{sn} - \rho_{iln} f_{sl})
+ \rho_{iln} \sum_{s\ne n,l} a_{sik}(\rho_{ikn} f_{sl} - \rho_{ikl} f_{sn})\cr
&+ \rho_{ikl} \sum_{s\ne k,l} a_{sin}(\rho_{iln} f_{sk} - \rho_{ikn} f_{sl})
+ \rho_{ikn}(\rho_{ikl}^2 \alpha_k - \rho_{iln}^2 \alpha_n)f_{kn}\cr
&+ \rho_{iln}(\rho_{ikn}^2 \alpha_n - \rho_{ikl}^2 \alpha_l)f_{ln}
+ \rho_{ikl}(\rho_{iln}^2 \alpha_l - \rho_{ikn}^2 \alpha_k)f_{lk}\cr
= &\ \rho_{ikn} \sum_{s\ne k,n} a_{sil}(\rho_{ikl} f_{sn} - \rho_{iln} f_{sl})
+ \rho_{iln} \sum_{s\ne n,l} a_{sik}(\rho_{ikn} f_{sl} - \rho_{ikl} f_{sn})\cr
&+ \rho_{ikl} \sum_{s\ne k,l} a_{sin}(\rho_{iln} f_{sk} - \rho_{ikn} f_{sl})\cr
&+ \alpha_n \rho_{ikn}\rho_{iln}(\rho_{ikn} f_{ln} - \rho_{iln} f_{kn})\cr
&+ \alpha_l \rho_{ikl}\rho_{iln}(\rho_{iln} f_{lk} - \rho_{ikl} f_{ln})\cr
&+ \alpha_k \rho_{ikl}\rho_{ikn}(\rho_{ikl} f_{kn} - \rho_{ikn} f_{lk})\cr}$$
This gives the desired lifting. \endpf
\bigskip
We now return to Step 1 in the proof of Theorem 3.1.  The lemma,
combined with Observation 3.3, shows that any syzygy on the elements
of ${\cal G}$ modulo $\langle x_5, x_0\rangle$ that can be expressed in
terms of the syzygies (3) can be lifted to a syzygy on ${\cal G}$, and
hence that no new elements of the Groebner basis will be produced
in those cases.  In fact,
{\it all} syzygies on ${\cal G}$ modulo $\langle x_5, x_0\rangle$ can
be expressed in terms of syzygies of the form (3), so
no new Groebner basis elements at all
are introduced after $G$.  For example, reducing $S(F_1', F_3')$
yields
$$x_4 F_1' - \rho_{134} x_2 F_3' - \rho_{124} x_4 F_4' \equiv 0
\hbox{ mod } \langle x_5, x_0 \rangle\ ,$$
or combining the two terms with $x_4$,
$$\rho_{134}(x_4 (\rho_{234} f_{12} - \rho_{124} f_{23}) -
x_2 (\rho_{234} f_{14} - \rho_{124} f_{34}))\equiv 0 \hbox{ mod }
\langle x_5, x_0 \rangle\leqno(6)$$
This relation is apparently of a different form than the ones in the Lemma,
but the fact that it too lifts it can be deduced from the Lemma as follows.
Modulo $\langle x_5, x_0 \rangle$, the first term on the
left of (6) is the ${\rho_{134}\over \rho_{123}}$ times
the term with the factor of $x_4$ on the left of relation (3) with $i = 2$.
Similarly, the second term on the left, modulo $\langle x_5, x_0 \rangle$,
is exactly the $x_2$ term on the left of relation (3) with $i = 4$.  We
form the corresponding linear combination of those two relations of the form
(3) and clear
denominators of $\rho_{123}$ yielding that
$$\eqalign{\rho_{123} \rho_{134}(x_4 (\rho_{234} f_{12}  - \rho_{124} f_{23})
&- x_2 (\rho_{234} f_{14} - \rho_{124} f_{34}))\cr
- \rho_{234} \rho_{124} (x_1 (\rho_{123} f_{34} - \rho_{134} f_{23}) &-
x_3(\rho_{134} f_{12} - \rho_{123} f_{14}))\cr}$$
equals the difference between the corresponding combination of the
right hand sides of the two relations (3).
But this shows that (6) lifts to a syzygy on the
quadrics in ${\cal G}$.  The other $S$-pairs are handled similarly
using relations (3) directly, and relations of this last type.

Our conclusion is that in the case that all $\rho_{ijk} \ne 0$, ${\cal G}$
is the reduced Groebner basis of $J$.  Computing the Hilbert function of
$J$ from this information, we see that $S = V(J)$ has degree 5 and codimension
3.  Indeed, $V(J) \cap V(x_0,x_5)$ consists of the five points with homogeneous
coordinates
$$
(1,0,0,0,0,0),\ (0,1,0,0,0,0),\ (0,0,1,0,0,0),\ (0,0,0,1,0,0),$$
$$\left(
{1\over \rho_{234}},\
{1\over \rho_{134}},\
{1\over \rho_{124}},\
{1\over \rho_{123}}, 0, 0
\right)
$$
so $S$ is reduced.

The remaining cases to consider are those where one or both of
$\rho_{134}, \rho_{234}$ are zero.  The arguments in those cases are
basically similar to the ones given here, so we will omit most of the
details and give only the form of the corresponding Groebner basis in each
case.  If {\it both} $\rho_{134} = \rho_{234} = 0$, then the quadrics in
(2) reduce to $$\{f_{13}, f_{14}, f_{23}, f_{24}, f_{34}\}\ .$$  Using the
vanishing of the Petri coefficients, we see that the initial ideal
of $J$ has the form
$$M_2 = \langle x_1 x_3, x_1 x_4, x_2 x_3, x_2 x_4, x_3 x_4, x_1^2 x_5\rangle$$
in this case.  Computing the Hilbert function gives degree 5 and codimension
3.

Finally if just one of the Petri coefficients,
say $\rho_{234}$, is zero, then from (2), we see that $J$ is generated
by
$$\rho_{134} f_{12} - \rho_{123} f_{14}, \
\rho_{124} f_{13} - \rho_{123} f_{14},\  f_{23},\  f_{24},\  f_{34}$$
The initial ideal of $J$ has the form
$$M_1 = \langle x_1 x_2, x_1 x_3, x_2 x_3, x_2 x_4, x_3 x_4, x_1 x_4^2\rangle$$
and the Hilbert function is the same as in the other cases.

This completes Step 1 of the proof of Theorem 3.1.
\bigskip
\noindent
{\bf Step 2.}
\bigskip
We now want to show that the family ${\cal S}$ of all surfaces $S$ of degree 5
obtained in step 1 is irreducible.
Looking at the minimal free resolution of the coordinate ring of $S$ in each
case,
we have that $J$ is a Gorenstein ideal of codimension three, since the
Betti diagram (as in the ``betti'' command of the Macaulay system of Bayer
and Stillman) is
$$\matrix{1&-&-&-\cr
          -&5&5&-\cr
          -&-&-&1\cr}$$
That is, writing $R = k[x_1,\cdots, x_5,x_0]$, the
minimal resolution has the form
$$0 \to R(-5) \to R(-3)^5 \to R(-2)^5 \to R \to R/J \to 0\leqno(7)$$

We can
use the structure theorem for Gorenstein ideals of codimension 3 ([BE])
to give a uniform description of the ideal $J$ valid in all cases.  Namely,
every ideal $J$ that appears here is generated by the $4 \times 4$
Pfaffians of a $5 \times 5$ skew-symmetric rank-4 matrix of linear forms
$A = (a_{ij})$ (the ``middle matrix'' of the resolution (7) under a
suitable choice of basis for $J$).  Apart from the requirement that
$\hbox{rank}(A) = 4$, the entries in $A$ are arbitrary.
It follows that the surfaces $S$ obtained in Step 1 form one irreducible
family, of dimension 35.  (This may also be seen by a normal bundle
calculation as in \S 2.)
\bigskip
\noindent
{\bf Step 3.}
\bigskip
\noindent
A general $C$ is contained in exactly one surface $S$ of the type
described above.
To complete the proof, we complete the basis (2) of $J$ to a basis of
$I(C)$.  Recall that we are assuming that $I(C)$ is generated by quadrics,
so any one further quadric in $I(C)$ not in $J$ (such as $f_{34}$ in
the case that all $\rho_{ijk} \ne 0$) will do the job.  Hence $C$ is
the complete intersection of $S$ and a quadric hypersurface.  (Conversely,
given a surface $S$ of degree 5 of the form above and a general quadric $Q$
-- not containing any component of $S$ for instance -- then $Q \cap S$ will be
a non-degenerate Gorenstein curve $C$ of degree 10 and arithmetic
genus 6 in ${\bf P}^5$.).  The family of surfaces $S$ is irreducible by Step 2,
and given $S$, to obtain $C \subset S$, the additional quadric $Q$ can be
chosen
essentially arbitrarily in $H^0(C, {\cal O}_S(2))$ which is irreducible of
dimension 15.
Hence, we have that ${{\cal H}_6}'$ is irreducible of dimension $35 + 15  =
50$. \endpf
\bigskip
(Note that, as expected, this is the same as the dimension of the family of all
smooth canonically-embedded curves in ${\bf P}^5$, which is
$$\hbox{dim}({\cal M}_6) + \hbox{dim}(PGL(6)) = 3\cdot 6 - 3 + 6^2 - 1 = 50\
.$$

We may ask if there is a result for $g = 6$ including the curves for which
$I(C)$ is not generated by quadrics, fully parallel to Proposition 2.1.  The
answer is yes, as we will now see.  However, the
situation is somewhat complicated by two new features.  First, even for
smooth curves, by the
classical Enriques-Petri theorem, there are two different possibilities for
the variety $V(f_{ij})$ when $I(C)$ is not generated by quadrics.  Indeed,
consider the family of ideals generated by quadrics in Petri's
form for which all $\rho_{ijk} = 0$ (so that $f_{ij}$ are a Groebner
basis for the ideal they generate), and in which, for simplicity,
the low-order terms in $f_{ij}$ are normalized to
$$q_{ij}(x_0, x_5) = b_{ij} x_0 x_5\ .\leqno(8)$$
Since $V(f_{ij},x_0,x_5)$ consists of the four points
$$
(1,0,0,0,0,0),\ (0,1,0,0,0,0),\ (0,0,1,0,0,0),\ (0,0,0,1,0,0),$$
we see that $V(f_{ij})$ is a reduced, arithmetically Cohen-Macaulay
surface of degree 4.  By the classification of surfaces of
degree $n-1$ in ${\bf P}^n$, we see that there are two components.
One, of dimension 11, has general point corresponding to the ideal of a
quartic scroll, another, of dimension 9 has general point corresponding to
the ideal of a Veronese surface.  (Recall that these Petri quadrics are
normalized so that the two additional points
$(0,0,0,0,1,0), (0,0,0,0,0,1)$ lie on the variety they define.
If we do not require these incidence conditions, the corresponding
components of the Hilbert scheme of surfaces have dimensions
$11 + 18 = 29$ (scrolls),  and $9 + 18 = 27$ (Veroneses) respectively.)

The curves lying on scrolls, and their degenerations,
and the plane quintics and their degenerations each form an irreducible
family, whose general element is a smooth curve.
This follows by an argument similar to that given above in Theorem 2.1.
For instance, on a scroll or a degeneration of a scroll, the canonical
curves are the residual intersections of the scroll and a cubic
containing 2 fixed lines.  This gives a
$$({5+3\choose 3} - 8 - 28) - 1 = 19$$
dimensional irreducible family of curves on each $S$.  Similarly,
on a Veronese surface or degeneration, the canonical curves belong
to a 20-dimensional irreducible family (e.g. the 2-uple images of
all the plane quintics on a smooth Veronese.)
 Hence,
an argument similar to the one given in the last
section of the proof of Theorem 1.2 shows that if $C$ is any
singular canonically-embedded curve on a scroll or a Veronese, or one of their
degenerations, then $[C]$ belongs to the same irreducible component of the
Hilbert scheme as the points in ${{\cal H}_6}'$.

There is a small second complication here as well.  Namely, by \S 3
of [S], in addition to the two cases we have accounted for,
in which the graded Betti number $\beta_{13} = 0$ ($C$ non-trigonal, non-plane
quintic), or $\beta_{13} = 3$ ($C$ trigonal, or plane quintic), there
is apparently another possibility when $g = 6$: $\beta_{13} = 1$.
By Proposition 3.2 of [S], this would happen only if there were exactly one
$\rho_{ijk} \ne 0$. The 4-vertex graph ${\cal T}$ introduced in Step 1
of the proof of Theorem 3.1 would be composed
of three edges forming a triangle together with one disconnected vertex.  We
follow the reasoning of Propositions 3.3 and 3.4 of [S].  Under
the assumption $\rho_{124} = \rho_{134} = \rho_{234} = 0$, but
$\rho_{123} \ne 0$, the quadrics $f_{12}, f_{13}, f_{23}$ depend
only on $x_1,x_2,x_3,x_0,x_5$ and satisfy the Petri syzygy $S_{123}$.
 From the other Petri syzygies, $a_{114} = a_{224} = a_{334}$, and
the linear form $x_4 - a_{114}$ must divide $f_{14}, f_{24}, f_{34}$.
We see from Theorem 1.4 of [S] that $f_{12}, f_{13}, f_{23}$,
together with the linear form $x_4 - a_{114}$ must
generate the ideal of a non-trigonal canonically-embedded curve of genus
5: $C_1 \subset H = V(x_4 - a_{114}) \cong {\bf P}^4 \subset {\bf P}^5$.
Furthermore,
$V = V(f_{ij}) = C_1 \cup P$, where $P$ is a plane.
By degree considerations, the
other component $C_2$ of $C$ lying in $P$ must be a conic in $P$.
(In order for $C$ to be Petri-general, $C_2$ must be reduced as well:
By the Petri construction, $(0,0,0,1,0,0)$ must be a smooth point of
$C$, but it cannot be
contained in $C_1$.  Hence it must be a smooth point of $C_2$.)

We claim that such a curve cannot be a canonically-embedded
curve (the condition ${\cal O}_C(1) \cong \omega_C$ will fail), even though
its Hilbert point almost certainly lies on the same component as
those of canonically-embedded
curves.  The reason is the following. The hyperplane $H$ containing
$C_1$ and the
plane $P$ containing $C_2$ must intersect transversely along a
line $L$ if $C$ is to span ${\bf P}^5$.  In order for $C$ to be
connected, we need $C_1 \cap C_2 \ne \emptyset$.

 By a standard fact
on dualizing sheaves ([C], Lemma 1.12),
$$\omega_C |_{C_2} \cong \omega_{C_2}\otimes
({\cal I}_{C_1}\otimes {\cal O}_{C_2})^{-1}\leqno(9)$$
where ${\cal I}_{C_1}$ is the ideal sheaf of $C_1$.  Since $C_1$ is
a non-trigonal curve, not containing $L$ as a component, $C_1 \cap L$
contains at most two points.  This implies that neither restriction
$\omega_C|{C_i}$ is correct.  For example,
in order for $C_2$ to embed as a conic by the dualizing
sheaf, $({\cal I}_{C_1}\otimes {\cal O}_{C_2})^{-1}$ would have
to have degree 4 on $C_2$, but this is impossible.

Hence, to obtain curves of this type, the only remaining
possibility is that the line $L$
is a component of $C_1$.  There are canonical curves
of genus 5 of this type: $C_1 = C_1' \cup C_1''$ where $C_1'$
has genus 3, $C_1'' \cong {\bf P}^1$, and $C_1'$ intersects $C_1''$
transversely
in three points $\{p_1,p_2,p_3\}$.  By (9), $\omega_C$ is very ample
in this case, and embeds $C_1'$ as a curve of degree 7 in ${\bf P}^4$,
with $C_1''$ as a trisecant line.  See 3.6 below for
a concrete example; if $p_1 + p_2 + p_3$ is not a divisor of
one of the $g_3^1$'s on $C_1'$, we even obtain curves of this type which
are not ``trigonal.''   However this case too leads to a situation
where $C$ fails to be canonically-embedded.  The reason is the
same as before.  In a general such curve, the smooth conic
$C_2$ will again meet the line $L$ in two points.  Hence
${\cal O}_{C_1''}(1)$ and ${\cal O}_{C_2}(1)$ are incorrect for a
canonically-embedded curve.
If $C_2$ meets $L$ at one or two of the points $p_i$, and $C_1'$ is smooth,
we obtain one or two non-planar triple points (with delta-invariant
$\delta = 2$).  These singularities are not even locally Gorenstein.
More degenerate curves also occur but the conclusion is the same in all cases.

In sum, the case $\beta_{13} = 1$ {\it does not} actually occur for
canonically-embedded curves when $g = 6$.  We have proved the following.
\bigskip
\proclaim Theorem 3.5.  Let ${\cal H}_6$ be the open subscheme
of the Hilbert scheme of curves of degree 10 and arithmetic
genus 6 in ${\bf P}^5$ corresponding to Petri-general canonically
embedded curves.   Then ${\cal H}_6$ is irreducible. \par

We conclude with two examples illustrating the analysis of the
$\beta_{13} = 1$ cases above.
\bigskip
\noindent
{\bf Example 3.6.}  Consider the curve $C = C_1 \cup C_2$, where $C_i$
are defined as follows.  Let $C_1$ be the non-trigonal genus 5 canonical curve
defined by the Petri-form quadrics
$$\eqalign{q_{12} &= x_1 x_2 - (x_5-x_0) x_1 + (x_0+x_5)x_3\cr
           q_{13} &= x_1 x_3 + (x_0+x_5)x_2\cr
           q_{23} &= x_2 x_3 + (x_0+x_5)x_1 + (4x_5-x_0)x_3\cr}$$
together with $x_4 = 0$.  (We have $\rho_{123} = 1$, and
$\alpha_1 = \alpha_2 = \alpha_3 = -(x_0 + x_5)$.  $C_1$ contains
$L = C_1'' = V(x_1, x_2, x_3, x_4)$ as a component. The other component
$C_1'$ is a curve of arithmetic genus 3 as above.

Next, let $C_2$ be the conic $V(x_1,x_2,x_3, x_4 x_5 - 2x_4 x_0 + 4 x_0 x_5)$.
Intersecting $I(C_1)$ and $I(C_2)$, we find a Groebner basis for
$I(C)$ yielding the following information.  The initial ideal of $I(C)$ is
$$\langle x_1 x_2, x_1 x_3, x_1 x_4, x_2 x_3, x_2 x_4, x_3 x_4,
x_1^2 x_5, x_2^2 x_5, x_4^2 x_5, x_3^3 x_5 \rangle$$
(note the {\it differences\/}
between this and the initial ideal given in Theorem 1.4 of [S]).
However, $C$ still has degree 10 and arithmetic genus 6.  The
Betti diagram of the minimal resolution of $I(C)$ has the form:
$$\matrix{1&-&-&-&-\cr
        -&6&6&1&-\cr
        -&1&6&6&1\cr
        -&-&-&1&1\cr}$$
We see that $\beta_{13} = 1$, but that
the resolution is not self-dual, which confirms the fact that $C$ is not
canonically-embedded.  However, we note that by a calculation,
$\hbox{dim }H^0(N_{C|{\bf P}^5}) = 50$, which
strongly suggests that the Hilbert point of $C$ lies on the
same component of the Hilbert scheme as the canonical curves of
genus 6.
\bigskip
In the case that $L$ is {\it not\/} a component of $C$, we obtain
non-canonically-embedded curves of degree 10 and arithmetic genus 6 for which
$\beta_{13} = 0$.
\bigskip
\noindent
{\bf Example 3.7.}  Consider the curve $C = C_1 \cup C_2$, where $C_i$
are defined as follows.  Let $C_1$ be the non-trigonal genus 5 canonical curve
defined by the Petri-form quadrics
$$\eqalign{q_{12} &= x_1 x_2 + (x_0+x_5)x_3 + x_0x_5\cr
           q_{13} &= x_1 x_3 + (x_0+x_5)x_2\cr
           q_{23} &= x_2 x_3 + (x_0+x_5)x_1\cr}$$
together with $x_4 = 0$.  (We have $\rho_{123} = 1$, and
$\alpha_1 = \alpha_2 = \alpha_3 = -(x_0 + x_5)$.  $C_1$ meets
$L = V(x_1, x_2, x_3, x_4)$ in the two points $(0,0,0,0,1,0), (0,0,0,0,0,1)$.)
Next, let $C_2$ be the conic $V(x_1,x_2,x_3, x_4 x_5 - 2x_4 x_0 + 4 x_0 x_5)$,
which meets $L$ in the same two points as $C_1$.
Intersecting $I(C_1)$ and $I(C_2)$, we find that instead of the quadric
$q_{12}$ above, $I(C)$ contains the quadric
$$f_{12} =  x_1 x_2 + (x_0+x_5)x_3 + (1/4 x_5-1/2 x_0)x_4 + x_0x_5$$
This shows that we have a curve similar to the
curves with $\rho_{123}, \rho_{124} \ne 0$
studied above in the proof of Theorem 3.1.  (Indeed the reader
will have no difficulty constructing a surface of degree 5 (a kind of
``degenerate
Del Pezzo surface'') containing $C$.  The initial ideal of $I(C)$  is
the same as in Example 3.6.  However, now the Betti diagram for the
minimal resolution of $I(C)$ has the form:
$$\matrix{1&-&-&-&-\cr
         -&6&5&1&-\cr
         -&-&6&6&1\cr
         -&-&-&1&1\cr}$$
We see that $\beta_{13} = 0$, but that as in the previous example,
the resolution is not self-dual. By another calculation,
$\hbox{dim }H^0(N_{C|{\bf P}^5}) = 50$, which again
strongly suggests that the Hilbert point of $C$ lies on the
same component of the Hilbert scheme as the canonical curves of
genus 6.

\vfill\eject
\noindent
{\bf References}
\bigskip
\frenchspacing
\item{[AH]} Arbarello, E. and Harris, J.  Canonical curves and quadrics
of rank 4, {\sl Compositio Math.} {\bf 43} (1981), 145-179.
\smallskip
\item{[ACGH]} Arbarello, E., Cornalba, M., Griffiths, P., and Harris, J.
{\sl Geometry of Algebraic Curves, Vol. 1}, New York: Springer Verlag, 1985.
\smallskip
\item{[BE]} Buchsbaum, D. and Eisenbud, D. Algebra structures for finite
free resolutions and some structure theorems for ideals of codimension 3,
{\sl Amer. J. of Math.} {\bf 99} (1977), 447-485.
\smallskip
\item{[C]} Catanese, F. Plurcanonical Gorenstein curves, {\it in} {\sl
Enumerative Geometry and Classical Algebraic Geometry}, {\sl Progress in
Mathematics}, {\bf 24}, Boston: Birkh\"auser, 1982.
\smallskip
\item{[KS]} Kollar, J. and Schreyer, F.-O. The Moduli of Curves is
Stably Rational for $g \le 6$, {\sl Duke Math. J.} {\bf 51} (1984), 239-242.
\smallskip
\item{[HM]} Herzog, J. and Miller, M.  Gorenstein Ideals of Deviation Two,
{\sl Communications in Algebra} {\bf 13} (1985), 1977-1990.
\smallskip
\item{[M]}  Mumford, D. {\sl Curves and their Jacobians}  Ann Arbor: University
of Michigan Press, 1975.
\smallskip
\item{[P]} Petri, K. \"Uber die invariante Darstellung Algebraischer Funktionen
einer Ver\"anderlichen, {\sl Math. Ann.} {\bf 88} (1923), 242-298.
\smallskip
\item{[SD]} Saint-Donat, B. On Petri's analysis of the linear system of
quadrics through a canonical curve, {\sl Math. Ann.} {\bf 206} (1973), 157-175.
\smallskip
\item{[S]} Schreyer, F.-O. A standard basis approach to syzygies of
canonical curves, {\sl J.f.d.reine u. angew. Math.} {\bf 421} (1991),
83-123.
\bigskip
\leftline{Department of Mathematics}
\leftline{College of the Holy Cross}
\leftline{Worcester, MA 01610}
\bigskip
\leftline{Current address, through 5/94:}
\bigskip
\leftline{Cornell University}
\leftline{Mathematical Sciences Institute}
\leftline{409 College Ave.}
\leftline{Ithaca, NY 14850}
\leftline{email: jlittle@msiadmin.cit.cornell.edu}

\bye